\documentclass{elsarticle}

\title{Enumerating combinatorial resultant decompositions of 2-connected rigidity circuits} 

\author[1]{Goran {Mali\'c}\fnref{fn1}}
\author[1]{Ileana Streinu\corref{cor1}\fnref{fn1}}

\cortext[cor1]{Corresponding author}
\fntext[fn1]{Funding: Both authors acknowledge funding from the NSF CCF:1703765 and CCF:2212309 grants to Ileana
	Streinu.}

\affiliation[1]{organization={Department of Computer Science, Smith College}, city={Northampton, MA},country={USA}}

\usepackage{hyperref}
\usepackage{xcolor}
\usepackage{amsfonts,amsmath,amsthm}
\usepackage{graphicx,tabularx}
\usepackage{epstopdf}
\usepackage{algpseudocode}
\usepackage{amssymb}
\usepackage{amstext}
\usepackage{tabularx}
\usepackage{multirow}
\usepackage[linesnumbered,ruled,vlined]{algorithm2e}\SetArgSty{}\SetFuncArgSty{} 
\usepackage{xspace}
\usepackage{caption}
\usepackage[capitalise]{cleveref}
\usepackage{subcaption}
\graphicspath{{./figures/}}

\newcounter{questioncounter}
\newtheorem{problem}[questioncounter]{Open Problem}

\newtheorem{lemma}{Lemma}
\newtheorem{theorem}{Theorem}

\usepackage{mdframed}
\newtheorem*{mainproblem}{Main Problem}

\newcommand{\cm}[1]{\operatorname{CM}_{#1}}

\newcommand{\grobner}{Gr\"obner}

\newcommand{\res}[3]{\operatorname{Res}(#1, #2, #3)}

\newcommand{\cres}[3]{\operatorname{CRes}(#1, #2, #3)}

\theoremstyle{definition}
\newtheorem{definition}{Definition}

\newcommand{\crt}{CR-tree\xspace}
\newcommand{\crts}{CR-trees\xspace}
\newcommand{\crd}{CR-decomposition\xspace}
\newcommand{\crds}{CR-decompositions\xspace}

\usepackage[margin=2.5cm]{geometry}

\begin{document}

\begin{abstract}
A rigidity circuit (in 2D) is a minimal dependent set in the rigidity matroid (in 2D), i.e.\ a minimal graph supporting a non-trivial stress in any generic placement of its vertices in $\mathbb R^2$.  Any rigidity circuit on $n\geq 5$ vertices can be obtained from rigidity circuits on a fewer number $4\leq n'<n$ of vertices by applying the combinatorial resultant (CR) operation \cite{malic:streinu:socg, malic:streinu:siaga}. The inverse operation is called a combinatorial resultant decomposition (CR-decomposition). Therefore, any rigidity circuit on $n\geq 5$ vertices can be successively decomposed into smaller circuits, until the complete graphs $K_4$ are reached. This sequence of CR-decompositions has the structure of a rooted binary tree called the combinatorial resultant tree (CR-tree).

A CR-tree encodes an \emph{elimination strategy} for computing circuit polynomials via Sylvester resultants. Different CR-trees lead to elimination strategies that can vary greatly in time and memory consumption. It is an open problem to establish criteria for optimal CR-trees, or at least to characterize those CR-trees that lead to good elimination strategies.

In \cite{malic:streinu:cgta} we presented an algorithm for enumerating CR-trees where we describe in detail the algorithms for decomposing 3-connected rigidity circuits in polynomial time. In this paper we focus on the rigidity circuits that are not 3-connected, which we simply call 2-connected.

In order to enumerate CR-decompositions of a 2-connected rigidity circuit $G$, a brute force exponential time search has to be performed among the subgraphs induced by the subsets of $V(G)$. This exponential time bottleneck is not present in the 3-connected case. In this paper we will argue that we do not have to account for all possible CR-decompositions of 2-connected rigidity circuits to find a good elimination strategy; we only have to account for those CR-decompositions that are a 2-split, all of which can be enumerated in polynomial time. We present algorithms and computational evidence in support of this heuristic.
\end{abstract}

	
\maketitle


\section{Introduction}

\noindent The \emph{combinatorial resultant} of two graphs $G_1=(V_1,E_1)$ and $G_2=(V_2,E_2)$ with respect to some $e\in E_1\cap E_2$ is the graph defined as $\cres{G_1}{G_2}{e}=(V_1\cup V_2, (E_1\cap E_2)-\{e\})$. If for a graph $G$ there exist graphs $G_1$ and $G_2$ such that $G=\cres{G_1}{G_2}{e}$, we say that $G$ has a \emph{combinatorial resultant decomposition} $(G_1,G_2,e)$, or \crd for short. In the context of rigidity theory, we have proven in \cite{malic:streinu:socg,malic:streinu:siaga} that any \emph{rigidity circuit}\footnote{All necessary definitions appear in \cref{sec:prelimRigidity}} $G$ with at least 5 vertices has a \crd $(G_1,G_2,e)$ such that both $G_1$ and $G_2$ are \emph{rigidity circuits}. The existence we have proven constructively by providing an algorithm: if $G$ is 3-connected, we obtain $G_1$ and $G_2$ by deleting certain vertices and adding an edge. If it is not 3-connected, we perform a \emph{2-split}, i.e.\ we find a separating pair $\{u,v\}$ and split $G$ along it into two components, and then adjoin a new edge $e=\{u,v\}$ to both. In other literature the edge $e$ is called a \emph{virtual edge} \cite{GutwengerMutzel, hopcroft:tarjan:73, tutte:connectivity:1966}; we refer to it as an \emph{elimination edge}.

Enumerating \crds of a rigidity circuit $G$ is a crucial step in our quest for computing \emph{circuit polynomials}, i.e.\ the unique, multi-variate and irreducible polynomials that completely encode the generic placements of rigidity circuits. Our choice of a name for the combinatorial resultant operation was not arbitrary. Combinatorial resultants provide an \emph{elimination strategy} for computing circuit polynomials: if $G=\cres{G_1}{G_2}{e}$ and both $G_1$ and $G_2$ are rigidity circuits whose circuit polynomials $p_{G_1}$ and $p_{G_2}$ are known, we have shown in \cite{malic:streinu:socg,malic:streinu:siaga} that we can compute the circuit polynomial for $G$ by computing the classical \emph{Sylvester resultant} of $p_{G_1}$ and $p_{G_2}$ with respect to the \emph{elimination variable} $x_e$ corresponding to the edge $e$ of $G_1$ and $G_2$.

Each rigidity circuit $G$ with at least 5 vertices can successively be decomposed into smaller and smaller rigidity circuits, until the smallest rigidity circuit, a $K_4$, is reached. This sequence of decompositions has the structure of a rooted binary tree that we call a \emph{combinatorial resultant tree} for $G$, or a \crt for short. \crts are not unique, each provides a different elimination strategy for computing the circuit polynomial for $G$ by building up from the circuit polynomial for a $K_4$ to increasingly complicated polynomials. How well a particular elimination strategy performs is very much dependent on its \crt. For a fixed circuit $G$ some \crts led to a fast computation of its circuit polynomial, while for other trees the computation would exhaust all resources and crash. It is therefore of interest to establish criteria for choosing a \crt that would result in an efficient elimination strategy.

The first natural step in this direction is to enumerate all \crts for a given rigidity circuit $G$. In \cite{malic:streinu:cgta} we presented an enumeration algorithm for \crts that builds the trees by decomposing $G$ in a certain way, depending on if $G$ is 3-connected or not. If it is 3-connected, then we can find all appropriate \crds $(G_1,G_2,e)$ of $G$ in polynomial time, where by appropriate we mean that both $G_1$ and $G_2$ have less vertices than $G$, and one of them is 3-connected with exactly $|V(G)|-1$ vertices. If $G$ is not 3-connected, then we have to perform a brute-force exponential time search among the subgraphs induced by subsets of $V(G)$. This exponential-time bottleneck is far from ideal, so we have restricted the enumeration of \crds for circuits that are not 3-connected to 2-splits only. We did not argue in detail in \cite{malic:streinu:cgta} why we make this choice; in this paper we will fill this gap.

\medskip\medskip
\noindent\fbox{
\parbox{\textwidth}{The \textbf{main goal} of this paper is to argue that for rigidity circuits that are not 3-connected it is sufficient to consider only 2-splits. In other words, we will argue that the elimination strategies given by 2-splits are better than other strategies.}
}

\medskip
By a better elimination strategy we mean a strategy that will outperform others with regard to elapsed time and the difficulty of computation. We present limited, but revealing evidence that 2-splits outperform other possible elimination strategies.

\medskip
\noindent
{\bf Overview of the paper.} All relevant definitions from 2D rigidity theory are presented in \cref{sec:prelimRigidity}. In the same section we define the Sylvester resultant and circuit polynomials. In section \cref{sec:mainConcepts} we give the necessary definitions and theorems for combinatorial resultants, \crds and \crts from \cite{malic:streinu:siaga}. From \cref{sec:2connected} until the end of the paper we work exclusively with rigidity circuits that are not 3-connected. We first give two algorithms for enumerating \crds of such circuits, a naive approach and an approach that considers 2-splits only. In \cref{sec:experiments} we present experimental evidence for our conjecture that \crds by 2-splits outperform other \crds. We conclude the paper with some final remarks in \cref{sec:concl}.

\section{Preliminaries: 2D rigidity circuits, the Sylvester resultant and circuit polynomials}
\label{sec:prelimRigidity}

We start by reviewing, in an informal way, those concepts and results from combinatorial rigidity theory in 2D that are relevant for our paper.

\medskip
\noindent
{\bf Notation.} Unless explicitly stated otherwise, the vertex set $V(G)$ of a graph $G$ is labeled with a subset of $[n]:=\{1,\dots,n\}$. The \emph{neighbours} $N(v)$ of $v$ are the vertices adjacent to $v$ in $G$. The support of $G$ is its set of edges $E(G)$. The \emph{vertex span} of a set of edges $E$ is the set of all vertices in $[n]$ incident to edges in $E$. A subgraph $H$ of $G$ is \emph{spanning} if its edge set $E(H)$ spans $[n]$. For $v\in V(G)$ we denote by $G-v$ the graph obtained by deleting $v$ and its incident edges from $G$. A \emph{non-edge} of $G$ is a pair $\{u,v\}\subset V(G)$ such that $\{u,v\}\notin E(G)$. For a non-edge $e$ of $G$ we define $G+e$ as the graph $(V(G),E(G)\cup e)$. We will often simplify notation by writing $G-v+e$ for the graph $(G-v)+e$.

\medskip
\noindent
{\bf Rigid graphs.} We say that a graph $G=(V,E)$ is \emph{rigid} if, up to rigid motions, the number possible placements with respect to any set of generic edge-lengths is finite, and \emph{flexible} otherwise. A \emph{minimally rigid} or {\em Laman} graph is a rigid graph in 2 dimensions such that the deletion of any edge results in a flexible graph. A graph with $n\geq 2$ vertices is a Laman graph if and only if it (a) has a total of $2n-3$ edges, and (b) it is (2,3)-sparse, i.e.\ no subset of $n'\leq n$ vertices spans more than $2n'-3$ edges. A \emph{Laman-plus-one graph} is obtained by adding one edge to a Laman graph: it has $2n-2$ edges and thus doesn't satisfy the $(2,3)$-sparsity condition (b).

\medskip
\noindent
{\bf Rigidity Matroids.} A matroid is a pair $(E,\mathcal I)$ of a finite ground set $E$ and a collection of subsets $\mathcal I$ of $E$ called \emph{independent sets}, which satisfy certain axioms \cite{Oxley:2011} which we skip (they are not relevant for our presentation). A \emph{base} is a maximal independent set. A set which is not independent is called \emph{dependent} and a \emph{minimal dependent set} is called a \emph{circuit}. For our purposes the following facts are relevant: (a) (hereditary property) a subset of an independent set is also independent; (b) all bases have the same cardinality, called the \emph{rank} of the matroid; (c) (the unique circuit property) if $I$ is independent and $I\cup \{e\}$ dependent for some $e\in E$, then there exists a unique circuit in $I\cup \{e\}$ that contains $e$.

\begin{figure}[ht]
	\centering
	\begin{tabular}{ccc}
		\includegraphics[width=.17\textwidth]{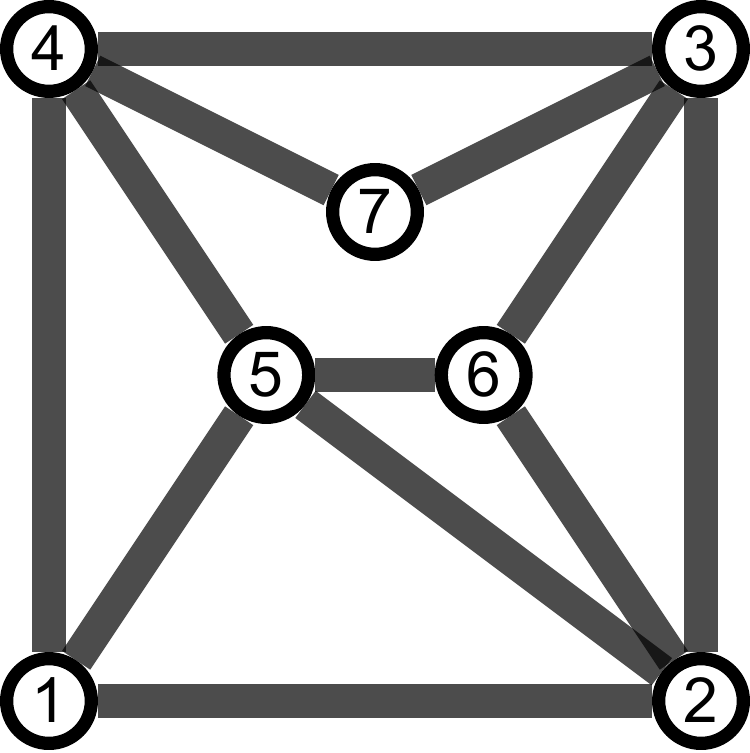} & \hspace{1.2cm} &  \includegraphics[width=.17\textwidth]{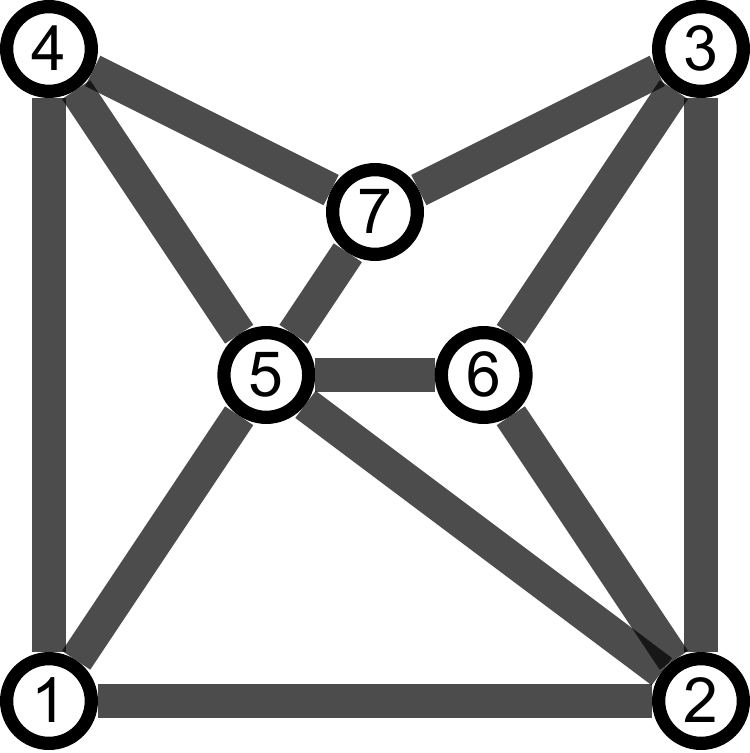}
	\end{tabular}
	\caption{Left: a Laman-plus-one graph on 7 vertices whose unique circuit is the Desargues-plus-one graph from bottom-right of \cref{fig:2sum}. This circuit is not spanning. Right: a circuit on 7 vertices. This circuit is spanning.}
	\label{fig:6circuits}
\end{figure}

The \emph{rigidity matroid} is defined on a ground set given by all the edges $E_n:=\{ij: 1\leq i < j \leq n\}$ of the complete graph $K_n$. Its bases are the Laman graphs on $n$ vertices. A \emph{Laman-plus-one} graph is a dependent graph obtained by adding an edge to a Laman graph, thus it has exactly $2n-2$ edges. A \emph{rigidity circuit} is a circuit in the rigidity matroid. It is a dependent graph with minimum edge support: deleting any edges leads to a Laman graph on its vertex set. In particular, a circuit in the rigidity matroid is a graph on $k\leq n$ vertices which (a) spans exactly $2k-2$ edges and (b) any subset of $k'<k$ vertices spans at most $2k'-3$ edges. A circuit whose vertex set is $V=[n]$ is said to be a \emph{spanning rigidity circuit} (Fig.~\ref{fig:6circuits}). A spanning rigidity circuit $C=(V,E)$ is a special case of a Laman-plus-one graph. Furthermore, a Laman-plus-one graph contains a unique circuit, which in general is not spanning.

In the rest of a paper rigidity circuits will simply be called circuits.

\medskip
\noindent
{\bf Connectivity.} Every circuit is necessarily $2$-connected, but may fail to be $3$-connected.  Let $G_1=(V_1,E_1)$ and $G_2=(V_2,E_2)$ be connected graphs, and denote by $V_{\cup}=V_1\cup V_2$ and $E_\cup=E_1\cup E_2$. If $G_1$ and $G_2$ are two graphs with exactly one common edge $uv$, we define their \emph{$2$-sum} as the graph $G=(V_{\cup},E_{\cup} \setminus \{uv\})$. The inverse operation of decomposing $G$ into $G_1$ and $G_2$ along the  non-edge $uv$ is called a $2$\emph{-split}.
The following Lemma will be important in the rest of the paper.

\begin{lemma}[\cite{BergJordan}]\label{lem:bergjordan}
	If $G_1$ and $G_2$ are rigidity circuits, then their 2-sum is a rigidity circuit. If $G$ is a circuit and $G_1$ and $G_2$ are obtained from $G$ via a 2-split, then $G_1$ and $G_2$ are rigidity circuits.
\end{lemma}

\begin{figure}[h]
	\centering
	\begin{tabular}[t]{ccc}
		\multirow{2}{*}[1.4cm]{\includegraphics[width=.17\textwidth]{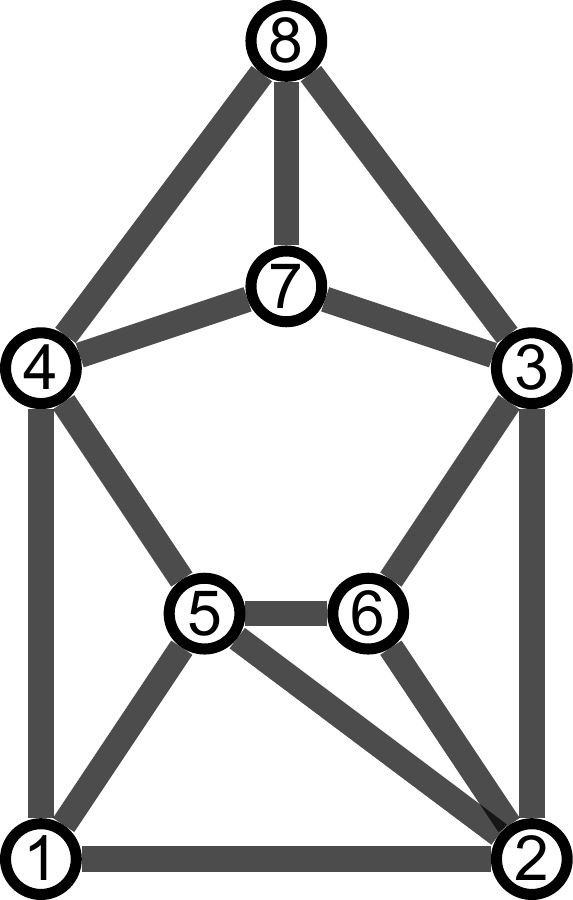}} & \hspace{1.2cm} & \includegraphics[width=.16\textwidth]{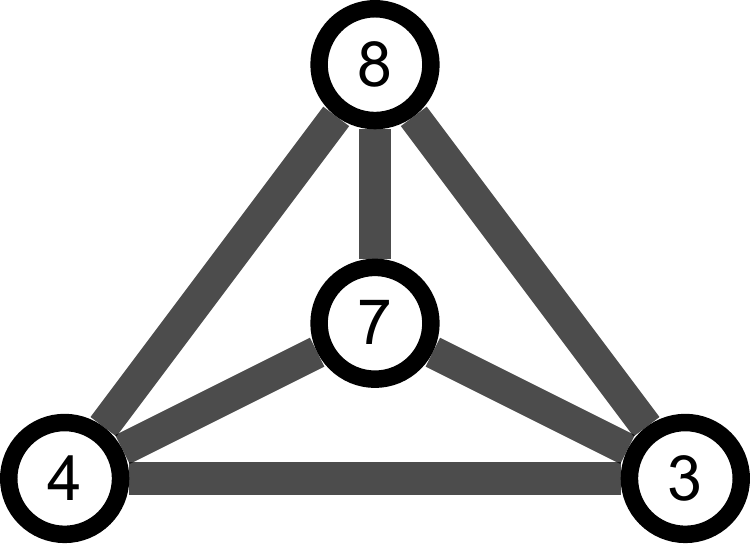}\\
		& \hspace{1.2cm} &\includegraphics[width=.17\textwidth]{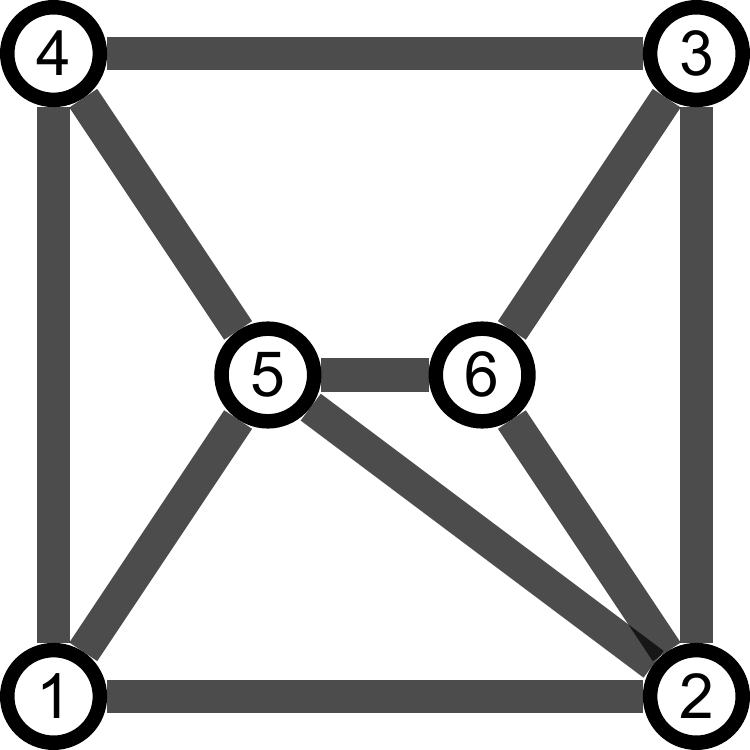}
	\end{tabular}
	\caption{Splitting a $2$-connected circuit (left) to get two $3$-connected circuits (right).}
	\label{fig:2sum}
\end{figure}

For an example, see \cref{fig:2sum}. A well-known theorem of Tutte \cite{tutte:connectivity:1966} implies that any circuit which is not $3$-connected can be uniquely decomposed into $3$-connected circuits with $2$-splits. A classical linear time algorithm for finding the 3-connected components was given by Hopcroft and Tarjan \cite{hopcroft:tarjan:73}.

\medskip
\noindent
\textbf{Circuit polynomials in the algebraic rigidity matroid.} Consider $n$ points in the Euclidean plane and let $x_{i,j}$ denote the squared distance between points $i$ and $j$. The \emph{Cayley-Menger matrix} on $n$ points is the $(n+1)\times (n+1)$ matrix obtained by augmenting the $n\times n$ Euclidean distance matrix on $n$ points with a top row $(0~1~\cdots~1)$ and a leftmost column $(0~1~\cdots~1)^\textnormal{T}$ of dimension $n+1$. The \emph{Cayley-Menger ideal} $\cm{n}$ on $n$ points is the prime ideal in $\mathbb Q[x_{i,j}\mid 1\leq i < j \leq n]$ generated by all $5\times 5$ minors of the Cayley-Menger matrix on $n$ points.

The algebraic rigidity matroid on $n$ points is a matroid on ground set $X_n=\{x_{i,j}\mid 1\leq i < j \leq n \}$. Its collection of independent sets $\mathcal I$ given by the subsets $I\subset X_n$ such that $\cm{n}\cap\mathbb Q[I] = \{0\}$. The algebraic rigdity matroid and the rigidity matroid are isomorphic. This result is folklore; we provided a proof in \cite{malic:streinu:siaga}.

A circuit in the algebraic rigidity matroid is a minimal $C\subset X_n$ such that $\cm{n}\cap\mathbb Q[C] \neq \{0\}$. It follows from a result of Dress and Lovasz \cite{DressLovasz} that for a circuit $C$ the ideal $\cm{n}\cap\mathbb Q[C]$ is  principal and generated by a polynomial $p_C\in \mathbb Q[C]$ unique up to multiplication with a rational number and irreducible over $\mathbb Q$. The polynomial $p_C$ is called the \emph{circuit polynomial} for $C$.

\medskip
\noindent\textbf{Sylvester resultant.} Let $f = a_0+a_1x+\cdots a_mx^m$ and $g = b_0+b_1x+\cdots b_nx^n$ be polynomials of degrees $m$ and $n$, with coefficients in a polynomial ring $R$. The Sylvester matrix $\operatorname{Syl}(f,g,x)$ of $f$ and $g$ is the $(m+n)\times(m+n)$ whose entries are the coefficients of $f$ and $g$:
\[
\operatorname{Syl}(f,g,x)=\begin{pmatrix}
	a_{m} & a_{m-1} & a_{m-2} & \cdots & a_0 & 0 & 0 & \cdots & 0\\
	0 & a_m & a_{m-1} & \cdots & a_1 & a_0 & 0 & \cdots & 0\\
	0 & 0 & a_m & \cdots & a_2 & a_1 & a_0 & \cdots & 0\\
	\vdots & \vdots & \vdots & \ddots & \vdots & \vdots & \vdots & \ddots & 0\\
	0 & 0 & 0 & \cdots & a_m & a_{m-1} & a_{m-2} & \cdots & a_0\\
	
	b_{n} & b_{n-1} & b_{n-2} & \cdots & b_0 & 0 & 0 & \cdots & 0\\
	0 & b_s & b_{n-1} & \cdots & b_1 & b_0 & 0 & \cdots & 0\\
	0 & 0 & b_s & \cdots & b_2 & b_1 & b_0 & \cdots & 0\\
	\vdots & \vdots & \vdots & \ddots & \vdots & \vdots & \vdots & \ddots & 0\\
	0 & 0 & 0 & \cdots & b_s & b_{n-1} & b_{n-2} & \cdots & b_0
\end{pmatrix}
\]
The submatrix consisting of the coefficients of $f$ is of dimension $n\times (m+n)$ and the submatrix consisting of the coefficients of $g$ is of dimension $m\times (m+n)$. Note that the columns align as show above only when $n=m$.

The \emph{Sylvester resultant}, denoted by $\res{f}{g}{x}$, is the determinant of the Sylvester matrix $\operatorname{Syl}(f,g,x)$. It is a polynomial in the ring $R$ given by a polynomial combination $a f + b g$ with $a,b\in R[x]$. We say that the Sylvester resultant \emph{eliminates} the variable $x$ from $f$ and $g$. For more details on the Sylvester resultant see \cite{CoxLittleOshea}.
We will be using the following lemma when discussing computational experiments.

\begin{lemma}\cite{malic:streinu:siaga}\label{lem:hom}
	Let $f,g\in R[x]$, where $R$ itself is a ring of polynomials. Denote by $\deg_x f$ and $\deg_x g$ the degrees in $x$ of $f$ and $g$, and let $h_f$ and $h_g$ denote their homogeneous degrees. The homogeneous degree of the Sylvester resultant $\res{f}{g}{x}$ is
	\[h_f\cdot\deg_xg + h_g\cdot\deg_xf -\deg_xf\cdot\deg_xg.\]
\end{lemma}

\section{Combinatorial resultants and combinatorial resultant trees}
\label{sec:mainConcepts}

This section reviews basic definitions on combinatorial resultants and combinatorial resultant trees from  \cite{malic:streinu:socg, malic:streinu:siaga}. Recall that $V_{\cup} = V_1\cup V_2$ and $E_{\cap}=E_1\cap E_2$.

\begin{definition}If $e\in E_1\cap E_2$ exists, we say that the graph $(V_{\cup}, E_{\cup}-\{e\})$ is the \emph{combinatorial resultant} of $G_1$ and $G_2$ on edge $e$, and denote by $\cres{G_1}{G_2}{e}$.\end{definition}

\noindent We will use the abbreviation CR in place of ``combinatorial resultant'' for the rest of the paper.

\medskip
\noindent
{\bf CRs of circuits.}
Given two circuits $G_1=(V_1,E_1)$ and $G_2=(V_2,E_2)$ with non-empty edge set intersection and $e\in E_1\cap E_2$, their CR with respect to edge $e$ is not always a circuit. When the common part $(V_\cap, E_\cap)$ is a Laman graph on its vertex set $V_\cap$, then the CR is necessarily a Laman-plus-one graph but not always a circuit.  Identifying the necessary and sufficient conditions for a CR of circuts to be a circuit is an open problem \cite{malic:streinu:socg, malic:streinu:siaga}. However, we will work exclusively with the case when the CR of two circuits is a circuit.

\medskip
\noindent
\begin{definition} If a circuit $G$ is obtained as a $G=\cres{G_1}{G_2}{e}$ of two other circuits $G_1$ and $G_2$, then we say that $(G_1,G_2,e)$ is a {\em combinatorial resultant decomposition} (\crd) of $G$. 
\end{definition}

\begin{theorem}[\cite{malic:streinu:socg, malic:streinu:siaga}]\label{thm:decomp}
	Every $3$-connected circuit $G$ with $n\geq 5$ vertices has a \crd $(G_1,G_2,e)$ such that $G_1$ is 3-connected and has $n-1$ vertices, and $G_2$ has at most $n-1$ vertices (but not necessarily 3-connected).
	
	A circuit that is not 3-connected has a \crd $(G_1,G_2,e)$ in which both $G_1$ and $G_2$ have at most $n-1$ vertices (and both might not be 3-connected). 
\end{theorem}

\noindent That $G_1$ is the circuit which is necessarily 3-connected when $G$ is 3-connected is merely a convention that we keep throughout this paper. In general we can interchange $G_1$ and $G_2$ arbitrarily.

\medskip
\noindent The proof of \cref{thm:decomp} for 3-connected circuits relies on the existence of non-adjacent pairs of \emph{admissible vertices} in 3-connected circuits $G$, i.e.\ vertices $v$ of degree 3 in $G$ such that for some non-edge $e$ of $G$ with endpoints in $N(v)$ the graph $G-v+e$  is again a 3-connected circuit \cite{BergJordan}. Such vertices may be absent from circuits that are not 3-connected (cf.\ \cref{sec:2connected}), in which case we establish the existence of a \crd via \cref{lem:bergjordan}.

\medskip
\noindent\cref{thm:decomp} provides an approach for computing circuit polynomials that avoids the costly \grobner{} basis computation with respect to an elimination order. In \cite{malic:streinu:socg,malic:streinu:siaga} we prove the following 
\begin{theorem}[\cite{malic:streinu:socg, malic:streinu:siaga}]\label{thm:res}
	Let $G$ be a circuit with a \crd $(G_1,G_2,e)$ and assume that the circuit polynomials $p_{G_1}$ and $p_{G_2}$ for $G_1$ and $G_2$ are known. Then the circuit polynomial $p_G$ for $G$ is an irreducible factor of $\res{p_{G_1}}{p_{G_2}}{x_e}$.
\end{theorem}

\noindent To determine which irreducible factor of $\res{p_{G_1}}{p_{G_2}}{x_e}$ is $p_G$ requires a factorization of the resultant and a test of membership in the Cayley-Menger ideal, the details of which are provided in \cite{malic:streinu:siaga}. \cref{thm:res} will be relevant in \cref{sec:experiments} where we discuss experimental results, however in the rest of this paper we will consider only the combinatorial aspects of \crds.

\medskip
\noindent
\textbf{Overview of the 3-connected \crd algorithm.} The proof of \cref{thm:decomp} is constructive, and a polynomial time algorithm for computing a \crd of a 3-connected circuit is presented in detail in \cite{malic:streinu:cgta}. For completeness we recall it here briefly: first identify a pair of non-adjacent admissible vertices $\{v,w\}$ of a circuit $G$. Set $G_1=G-v+e$ where $e$ is a non-edge of $G$ with vertices in $N(v)$ such that $G_1$ is a 3-connected circuit (at least one such $e$ is guaranteed to exist \cite{BergJordan}). For $G_2$ consider first $G - w + e$. The graph $G-w$ is a Laman graph, and therefore a basis in the rigidity matroid. Therefore $G-w+e$ is dependent and there exists a unique circuit in $G-w+e$ containing $e$. We take $G_2$ to be that unique circuit.

\medskip
\noindent
A \crd is not unique, and moreover not all \crds are captured  by this theorem; certain circuits $G$ can have a \crd $(G_1,G_2,e)$ such that one or both of $G_1$ or $G_2$ have the same number of vertices as $G$. However, we are interested in inductive constructions from which smaller circuits are used to build larger ones, hence for that purpose we are interested only in the \crds which produce circuits with fewer vertices than $G$.

\medskip
\noindent
{\bf\crts.}
Given a rigidity circuit $G$ with $n\geq 5$ we can construct a tree by repeatedly applying \cref{thm:decomp} to 3-connected circuits or 2-splits to 2-connected circuits. This construction stops once all leaves are $K_4$ circuits. We call this tree a \emph{\crt}.

\begin{definition}\label{def:resTree}
	Let $G$ be a rigidity circuit with $n\geq 5$ vertices. A \emph{combinatorial resultant tree (\crt)} for $G$ is a rooted binary tree with its root labeled with $G$, such that:
	\begin{enumerate}
		\item[\emph{(i)}] each node is labeled with a rigidity circuit;
		\item[\emph{(ii)}] leaves are labeled with graphs isomorphic to $K_4$;
		\item[\emph{(iii)}] children of non-leaf nodes $G'$ have labels $L$ and $R$ which satisfy $G'=\cres{L}{R}{e}$ for some $e\in L\cap R$ and both have fewer vertices than $G'$.	
	\end{enumerate}
\end{definition}

Recall that \cref{thm:decomp} states that for 3-connected circuits $G$ with a \crd $(G_1,G_2,e)$, the circuit $G_1$ is the one that has $n-1$ vertices and is 3-connected by convention. In this definition there is no requirement to keep that convention, or to even make a consistent choice that specifies which one of $L$ or $R$ is 3-connected with $n-1$ vertices. Furthermore, note that in \cite{malic:streinu:siaga}, this is referred to as a {\em combinatorial circuit resultant tree}.

A \crt is not unique, and already for circuits with a small number of vertices the number of different possibilities grows rapidly \cite{malic:streinu:cgta}.

\medskip
\noindent
\textbf{Algorithmic aspects and data structures.} We have presented an algorithm for enumerating the \crts of a 3-connected circuit in \cite{malic:streinu:cgta}. It relies on a memoization technique which computes the \crds only of those circuits that appear for the first time up to isomorphism. Once an isomorphic instance $H$ of a previously decomposed circuit $G$ is encountered, the \crt is truncated at $H$, i.e.\ $H$ becomes a leaf in this tree. Hence in \cite{malic:streinu:cgta} we are enumerating \emph{truncated CR-trees}.

\begin{definition}\label{def:truncResTree} 
	Let $G$ be a rigidity circuit with $n\geq 5$ vertices. The \emph{first truncated \crt} for $G$ is obtained as follows. Let $\mathcal C$ be the collection of pairwise non-isomorphic circuits, currently containing only $G$.
	\begin{enumerate}
		\item[\emph{(i)}] Compute a \crd $(L,R,e)$ for the root $G$. If $L$ and $R$ are isomorphic, mark one as truncated (e.g.\ $R$) and add the other (e.g.\ $L$) to $\mathcal C$. Otherwise add both to $C$. Proceed in a depth-first manner to the first child that hasn't been truncated.
		\item[\emph{(ii)}] For any circuit $G'$ computed so far that has not been marked as truncated, compute a \crd $(L',R',e')$. If $L'$ (resp.\ $R'$) is not isomorphic to any circuit in $\mathcal C$, add it to $\mathcal C$. Otherwise, mark it as truncated. Proceed to the first node that hasn't been truncated in a depth-first manner.
		\item[\emph{(iii)}] Repeat \emph{(ii)} in a depth-first manner until all leaves are a $K_4$ or marked as truncated.
	\end{enumerate}
	All remaining truncated \crts for $G$ are computed in the same way, with the exception that now in step \emph{(i)} both children could be isomorphic to circuits in $\mathcal C$.
\end{definition}

Truncated \crts are motivated as a method of reducing duplicate data. Ideally we want to enumerate isomorphism classes of CR-trees, where $f\colon T_1\to T_2$ is defined as a \crt isomorphism if it is an isomorphism of rooted binary trees and such that siblings in $T_1$ are isomorphic to siblings in $T_2$. More precisely, if $(G_1,G_2)$ are siblings in $T_1$ on nodes $N_1$ and $N_2$, then for siblings $(H_1,H_2)$ on nodes $f(N_1)$ and $f(N_2)$ either (i) $G_1\cong H_1$ and $G_2\cong H_2$, or (ii) $G_1\cong H_2$ and $G_2\cong H_1$, where $\cong$ denotes graph isomorhpism. However, computing all such isomorphisms is costly and truncated \crts mitigate that cost while still reducing the duplicate data to an extent.

Note that if we keep the convention in which $L$ is always 3-connected with one vertex fewer than its parent, then the left-most branch of the first truncated \crt is never truncated and therefore of maximal depth.

\medskip
\noindent
A truncated \crt is represented by a data structure that keeps a list of (pointers to) \emph{circuit nodes} (C-nodes), and a list of (pointers to) \emph{branching nodes} (B-nodes). A B-node encodes a single \crd of a circuit $G$ as a triple of pointers to $G$, its left child $L$ and its right child $R$. A C-node points to a unique circuit $G$ and to a list of B-nodes (one B-node for each possible \crd of $G$). Moreover it also stores an isomorphism of $G$ so that proper relabeling can be applied when necessary. Therefore the space of truncated \crts for a circuit $G$ is completely encoded by a list of unique circuits appearing in all possible \crts of $G$, a list of C-nodes (each C-node pointing to a unique circuit) and a list of B-nodes (one for each \crd of each unique circuit). This data is sufficient to reconstruct any \crt of $G$. Complete details are given in \cite{malic:streinu:cgta}.

\section{\bf The CR-decomposition of a 2-connected circuit}\label{sec:2connected}
For the rest of the paper we consider only those rigidity circuits $G$ that are 2- but not 3-connected. However, at times we will have to compare 2- but not 3-connected circuits to 3-connected circuits; hence, for the purpose of brevity, when we say that a circuit is 2-connected, we will also mean that it is \emph{not} 3-connected.

\begin{figure}[ht]
	\centering
	\begin{tabular}{ccc}
		\includegraphics[width=.35\textwidth]{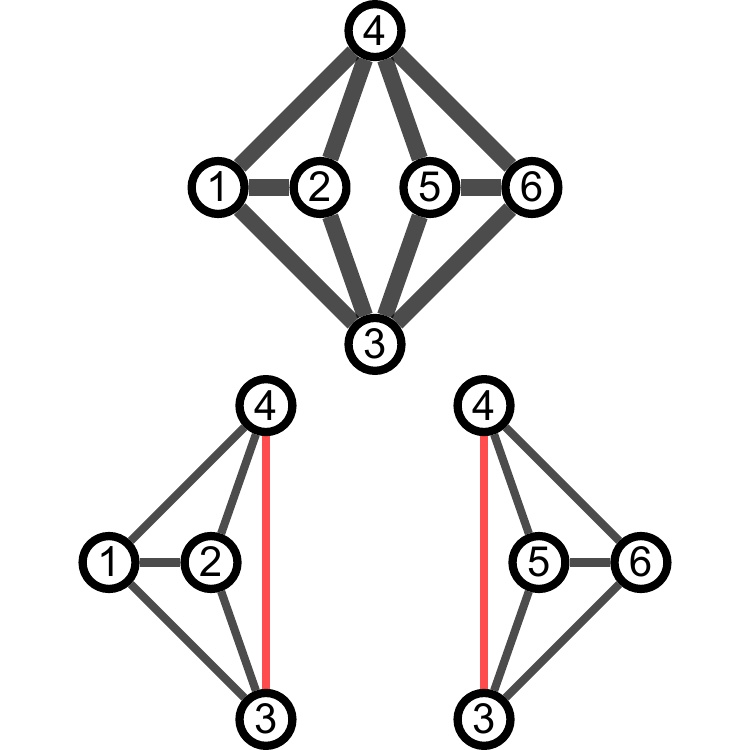} & \hspace{1cm} & \includegraphics[width=.35\textwidth]{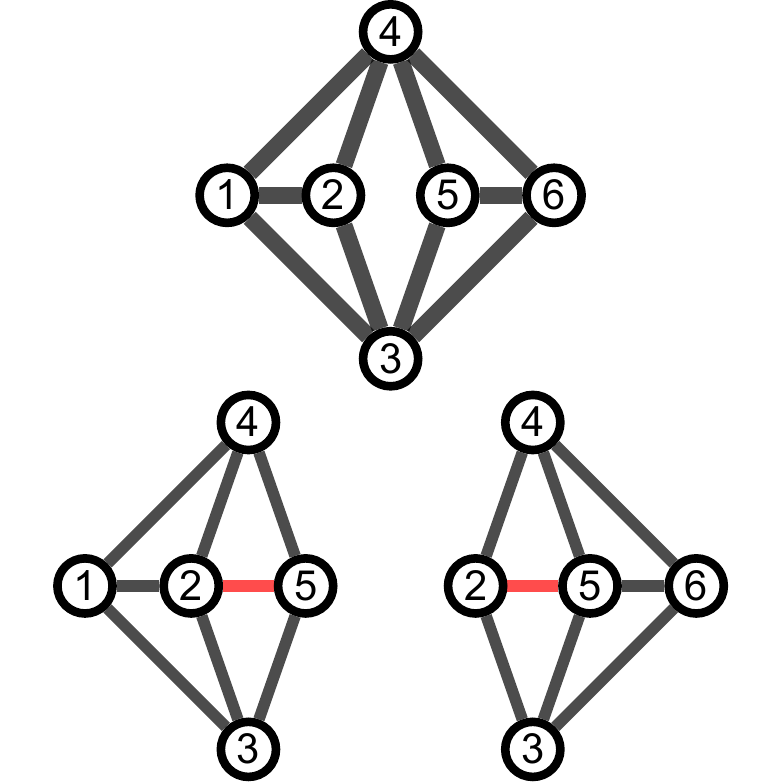}
	\end{tabular}
	\caption{Left: a \crd of the double banana into two $K_4$ graphs which is a 2-split. Right: a \crd of the double banana into two wheels on 4 vertices with the double triangle 2354 in common. The elimination edge for both cases is shown in red.}\label{fig:doubleBananaW4}
\end{figure}

Recall from the previous section that an admissible pair for $G$ is a vertex $v$ of degree 3 and a non-edge $e$ of $G$ with vertices in $N(v)$ such that $G-v+e$ is 3-connected. The main obstacle to applying the \crd algorithm to 2-connected circuits is that admissible pairs might not exist. This property is exhibited already on the smallest 2-connected circuit, the  ``\emph{double banana}'' shown in the top row of \cref{fig:doubleBananaW4}. However, it is still possible to decompose the double banana into smaller circuits in two ways: as a 2-sum of two $K_4$ circuits, or as a combinatorial resultant of two 4-wheels, as shown on the bottom of \cref{fig:doubleBananaW4}.

Recall that by \cref{lem:bergjordan} any 2-split of $G$ results in two rigidity circuits and by a theorem of Tutte \cite{tutte:connectivity:1966} $G$ can eventually be decomposed into unique 3-connected components via 2-splits. A 2-split is a CR-decomposition $(G_1,G_2,e)$ where $e$ is a non-edge between a separating pair of $G$, hence a \crt for $G$ can be obtained by first splitting $G$ into the 3-connected components given by Tutte's theorem, and then by decomposing each 3-connected component with the \crd algorithm until we reach $K_4$ graphs. However, $G$ can have \crds into 3-connected circuits that are not 2-splits, as shown in \cref{fig:doubleBananaW4}, which can lead to different \crts for $G$. For an efficient algorithm that can find all \crds of 2-connected circuits, we require necessary and sufficient criteria to establish when a CR $\cres{G_1}{G_2}{e}$ is a circuit. The graphs $(G_1,G_2,e)$ in a \crd of $G$ necessarily have to intersect on a Laman graph by following a sparsity count, however we are not aware of any sufficient criteria.

\begin{problem}\label{problem:2conn}
	Find necessary and sufficient conditions for a 2-connected rigidity circuit on $n\geq 6$ vertices to have a \crd other than a 2-split. 
\end{problem}

We therefore distinguish between two approaches, a naive exponential-time approach which brute-force enumerates all possible \crds of a 2-connected circuit in \cref{sec:naive}, and a partial, linear-time approach which enumerates only those \crds that are 2-splits in \cref{sec:split}.

\medskip

To motivate the approach which enumerates only the 2-splits, recall that the problem that motivates the enumeration of \crts is to find a good algebraic elimination strategy when computing circuit polynomials. We do not need to necessarily know all possible \crds in order to compute a circuit polynomial, we only need one strategy, or a small set of strategies that will in general outperform all others. Experimental results (c.f.\ \cref{sec:experiments}), although currently limited, suggest that decomposing a 2-connected circuit via a 2-split is the better strategy. At this point we can only conjecture that a \crd $(G_1,G_2,e)$ that is a 2-split will result in a ``simpler'' algebraic computation because the variables of the corresponding circuit polynomials are already separated. That is, the only common variable in the supports of the circuit polynomials $p_{G_1}$ and $p_{G_2}$ is $x_e$, whereas in a \crd that is not a 2-split, there might be a large number of variables in common.

\subsection{\textbf{The naive approach for decomposing a 2-connected circuit}}\label{sec:naive} The naive approach is illustrated in \cref{alg:naive}. It first enumerates all possible \crds $(G_1,G_2,e)$ of the input graph, and then checks if both $G_1$ and $G_2$ are circuits. This check is performed with the Decision Pebble Game of \cite{streinu:lee:pebbleGames:2008}.

\begin{algorithm}[ht]
	\caption{Naive enumeration of CR-decompositions of a 2-connected circuit}
	\label{alg:naive}
	\KwIn{A 2-connected circuit $G$ with $n \geq 5$ vertices}
	\KwOut{Set of all CR-decompositions for $G$}
	\SetKwFunction{pebbleGame}{PEBBLE-GAME}
	\BlankLine
	allDecompositions = empty\;
	\For{each $n'$-subset $V'$ of $V(G)$, $2\leq n'< n-1$}{
		$G' = $ subgraph spanned by $V'$\;
		\If{$|E(G')| = 2n'-3$}{
			\For{each non-trivial bipartition $A\cup B$ of $V\setminus V'$ and a non-edge $e$ of $G'$}
			{
				$G_1 = $ graph spanned by $V'\cup A$\;
				\textbf{if} \pebbleGame{$G_1+e$} == False \textbf{then} break\;
				$G_2 = $ graph spanned by $V'\cup B$\;
				\textbf{if} \pebbleGame{$G_2+e$} == False \textbf{then} break\;
				append $(G_1+e,G_2+e,e)$ to allDecompositions
			}
		}
	}
	\Return allDecompositions
\end{algorithm}

This approach starts by checking for each $n'$-subset $V'$ of vertices, where $2 \leq n' < n-1$, whether it spans a Laman subgraph $G'$. This step is required because the graphs $(G_1,G_2,e)$ in any \crd intersect on a Laman graph. Because circuits are (2,3)-sparse, checking whether $V'$ spans a Laman graphs amounts to checking that $n'=2n-3$.  
Then, every non-trivial bipartition $A\cup B$ of the remaining $n-n'$ vertices in $V\setminus V'$ gives rise to two graphs $G_1$ and $G_2$ spanned by the vertices in $V'\cup A$ and $V'\cup B$, respectively. Together with a non-edge $e$ of the common Laman subgraph $G'$, $(G_1,G_2,e)$ gives rise to a potential \crd. This triple will be a \crd if $G_1$ and $G_2$ are both circuits, which can be checked with the Decision Pebble Game algorithm.

\medskip
\noindent It is evident that the naive approach requires exponential time in order to check all possibilities for a common Laman subgraph. This is in contrast to the 3-connected case where the algorithm for enumerating all \crds runs in  $O(n^4)$ time. \cite{malic:streinu:cgta}.

\subsection{\textbf{Decomposing a circuit via 2-splits}}\label{sec:split}



\cref{alg:split} illustrates the procedure of enumerating all \crds that are a 2-split of the input circuit. We first enumerate all separating pairs of the input circuit, and then split it accordingly. The enumeration of separating pairs is performed by the linear time SPQR-tree algorithm of Gutwenger and Mutzel \cite{GutwengerMutzel} which in general can find all separating pairs of a multi-graph with $V$ vertices and $E$ edges in $O(|V|+|E|)$ time. Their algorithm is based on the SPQR-tree data structure introduced by Di Battista and Tamassia \cite{DiBattista89, DiBattista90}, originally in the context of planarity testing.

\begin{algorithm}[ht]
	\caption{Enumeration of CR-decompositions of a 2-connected circuit by 2-splits}
	\label{alg:split}
	\KwIn{A 2-connected circuit $G$ with $n \geq 5$ vertices}
	\KwOut{Set of all CR-decompositions for $G$ that are 2-splits}
	\BlankLine
	allDecompositions = empty\;
	SP = all separating pairs computed by the SPQR-tree algorithm \cite{GutwengerMutzel}\;
	\For{each pair $\{u,v\}\in\textnormal{SP}$}
	{
			$(G_1,G_2) = $ the two splits of $G$ with respect to $\{u,v\}$\;
			$e = \{u,v\}$\;
			append $(G_1,G_2,e)$ to allDecompositions
	}
	\Return allDecompositions
\end{algorithm}

Note that we do not have to test if $G_1$ and $G_2$ are circuts, as was the case in the naive approach. This is ensured by \cref{lem:bergjordan}.

\begin{figure}[ht]
	\centering
	\includegraphics[trim={0 10.5cm 0 2cm},width=.5\textwidth]{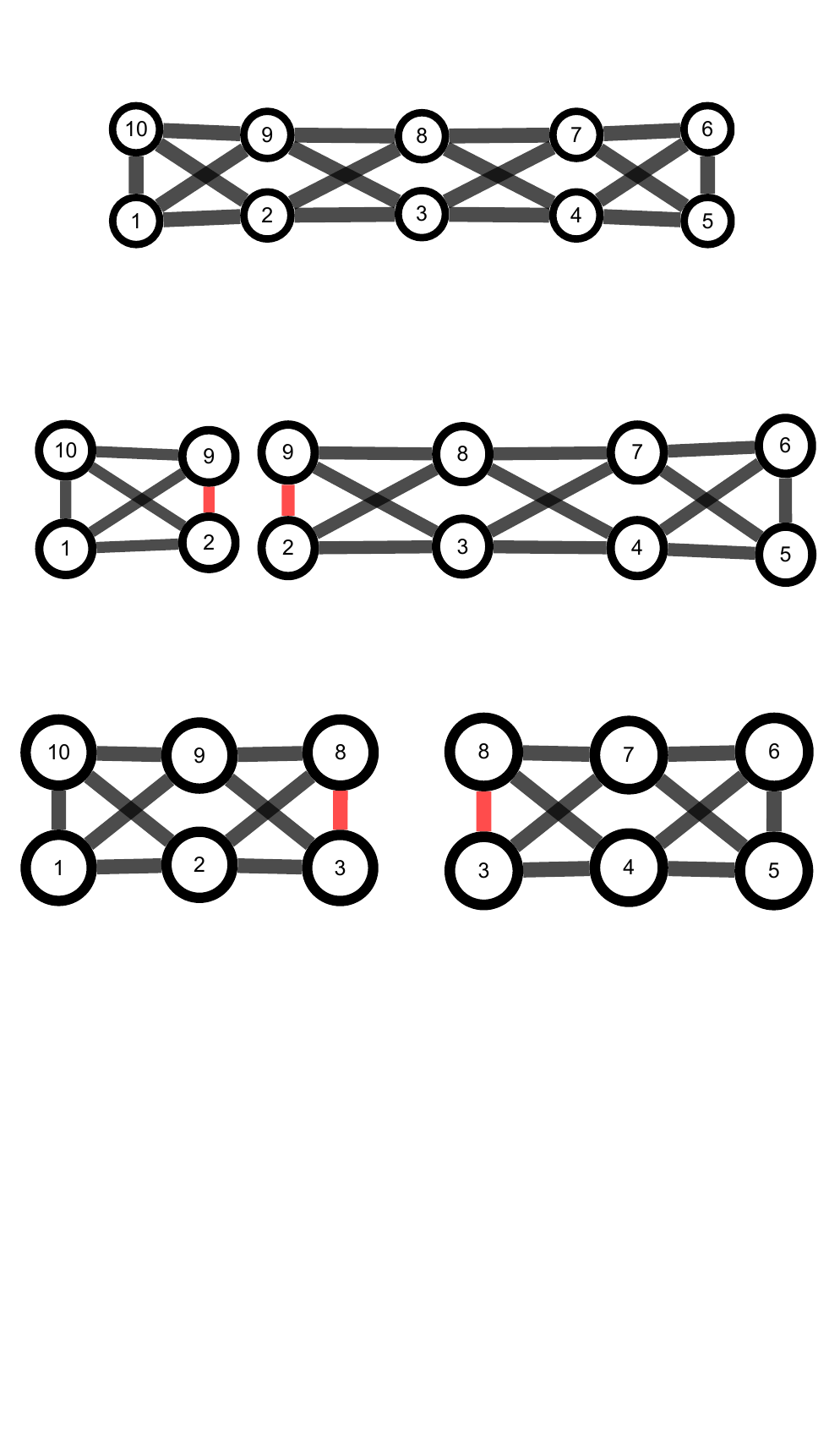}
	\caption{From top to bottom: a 2-connected circuit on 10 vertices with 3 separating pairs $\{2,9\}$, $\{3,8\}$ and $\{4,7\}$; a 2-split along the separating pair $\{2,9\}$; a 2-split along the separating pair $\{3,8\}$.}\label{fig:10v}
\end{figure}

\medskip
\noindent\textbf{Example and heuristics.} Consider the 2-connected circuit shown in the top row of \cref{fig:10v}.
This circuit has 3 separating pairs $\{2,9\}$, $\{3,8\}$ and $\{4,7\}$. The pairs $\{2,9\}$ and $\{4,7\}$ give rise to two isomorphic decompositions, so we consider only the \crds with respect to the non-edges $\{2,9\}$ and $\{3,8\}$.

\emph{Case $\{2,9\}$.} This case is shown on the left of \cref{fig:10v}. The left circuit is a $K_4$, and the right circuit is again 2-connected. We have computed the circuit polynomials for both: the circuit polynomial for $K_4$ has homogeneous degree 3 and is of degree 2 in each variable, and for the right circuit the circuit polynomial has homogeneous degree 20 and is of degree 8 in each variable \cite{malic:streinu:GitHubRepo}. This data gives us 
\begin{enumerate}
	\item[(i)] the dimension of the Sylvester determinant, which is given by the sum of the degrees of the elimination variable in the input polynomials. In this case the dimension is $(8+2)\times(8+2) = 10\times 10$ for any common variable;
	\item[(ii)] an upper bound for the homogenous degree of the circuit polynomial, which is given by the homogeneous degree of the Sylvester resultant (\cref{lem:hom}). In this case the upper bound is 48.
\end{enumerate}

\emph{Case $\{3,8\}$.} This case is shown on the right of \cref{fig:10v}. Here both circuits in the \crd are 2-connected and isomorphic to the double banana. The circuit polynomial for the double banana has homogeneous degree 8 and is of degree 4 in each variable. Therefore
\begin{enumerate}
	\item[(i)] the dimension of the Sylvester determinant in this case is $(4+4)\times(4+4) = 8\times 8$.
	\item[(ii)] an upper bound for the homogeneous degree of the circuit polynomial is 48.
\end{enumerate}

Which is the better choice? The advantage of Case $\{2,9\}$ is that the circuit polynomial $p_{K_4}$ for $K_4$ is the simplest possible, it has only 22 monomial terms. However, the circuit polynomial for the other circuit has 9 223 437 monomial terms, potentially negating any advantage that $p_{K_4}$ might bring.

The advantage of Case $\{3,8\}$ is that the Sylvester determinant is of smaller dimension, $8\times 8$ compared to $10\times 10$ of Case $\{2,9\}$. Furthermore, the two circuit polynomials up to change of labels are the same polynomial with 1752 monomial terms. This information is relevant for optimization of space complexity, as we don't have to keep both polynomials in memory, but only one of them together with a permutation encoding the change of labels. 

Unfortunately we were not able to compute any of the two cases, however we are inclined to guess that Case $\{3,8\}$ is the better choice because the Sylvester determinant is of lower dimension and with relatively simpler entries. This determinant will have 40 non-zero entries, each being a polynomial with at most 869 monomial terms, whereas in Case $\{2,9\}$ the large circuit polynomial is contributing 18 non-zero entries, some of which have more than 3 million monomial terms.

This is in contrast to our heuristics for the 3-connected case in which we always prefer a decomposition with a $K_4$ over all others, if there is one. The reason for this is as follows: consider a 3-connected circuit $G$ with $n$ vertices and with two \crds $(G_1,K_4,e_1)$ and $(G_2,G_3,e_2)$ where neither $G_2$ or $G_3$ are isomorphic to a $K_4$. The circuit $G_1$ has exactly $n-1$ vertices and so does one of $G_2$ or $G_3$, say $G_2$. We expect $G_1$ and $G_2$ to have relatively similar circuit polynomials, but $G_3$ could have anything between 5 and $n-1$ vertices and would therefore be much larger than the circuit polynomial for the $K_4$ circuit. That this doesn't necessarily happen with 2-connected circuits is an interesting phenomenon warranting further experimentation.

\section{Experimental observations regarding 2-connected circuits.}\label{sec:experiments} We conjecture that a decomposition of $G$ by 2-splits is the \emph{best possible} in the sense that the corresponding circuit polynomials are the least complicated with respect to the number of indeterminates, number of monomial terms, degrees in individual indeterminates, or the homogeneous degree. We have limited data supporting this conjecture given in Examples 1--5.

Recall that to find a circuit polynomial via \cref{thm:res} we have to compute the Sylvester resultant and establish which of its irreducible factors is the circuit polynomial. In all examples provided the Sylvester resultants were all shown to be irreducible, therefore they are the sought circuit polynomials. In this paper we note only the (wall-clock) time it took to compute the Sylvester resultant. We omit the timing of the irreducibility test since in general one has to perform a further computation to extract the circuit polynomial, namely a factorization and an ideal membership test \cite{malic:streinu:socg,malic:streinu:siaga}.

In all examples the computational time refers to the wall-clock time measured by the RepeatedTiming function in Mathematica V13.3.1, i.e.\ it is the average of 5 trials with best and worst trial scrubbed. All computations were performed on a machine with the  AMD Ryzen 9 5950x 16-core CPU (with a base 3.4GHz clock overclocked to 5GHz) and 64GB of DDR4 RAM.

\medskip
\noindent\textbf{Example 1} \cref{tbl:dblBan} compares the computation of the circuit polynomial for the double banana as shown in \cref{fig:doubleBananaW4}. The computation corresponding to the 2-split was completed within 0.02 seconds, producing an irreducible polynomial of homogeneous degree 8 with 1762 monomial terms.

However, the computation corresponding to the \crd given by two wheels of 4 vertices with a double triangle 2354 in common crashed after memory capacity was exhausted. The homogeneous degree could be computed because we can apply the formula for the homogeneous degree of the Sylvester resultant determinant $\operatorname{Syl}(f,g,x)$ given by $h_fd_g + h_gd_f -d_fd_g$ where $h_f$ (resp.\ $h_g$) is the homogeneous degree of $f$ (resp.\ $g$) and $d_f$ (resp.\ $d_g$) is the degree in $x$ of $f$ (resp.\ $g$) \cite{malic:streinu:siaga}.

\begin{table}[h]
	\centering
	\begin{tabular}[t]{|l|c|}
		\hline
		\multicolumn{2}{|c|}{Elimination guided by the 2-split in \cref{fig:doubleBananaW4} Left}\\
		\hline
		Syl Det Dimension & $4\times 4$\\
		\hline
		Syl Resultant & Hom.\ poly.\ of hom.\ deg.\ 8 \\
		\hline
		\# of mon.\ terms & 1752 \\
		\hline
		Computational time & 0.017 seconds \\
		\hline
	\end{tabular}\vspace{.5cm}
	\begin{tabular}[t]{|l|c|}
		\hline
		\multicolumn{2}{|c|}{Elimination guided by the \crd in \cref{fig:doubleBananaW4} Right}\\
		\hline
		Syl Det Dimension & $8\times 8$\\
		\hline
		Syl Resultant & Hom.\ poly.\ of hom.\ deg.\ 48 \\
		\hline
		\# of mon.\ terms & Computation exceeded memory capacity \\
		\hline
		Computational time & Computation exceeded memory capacity \\
		\hline
	\end{tabular}
	\caption{Comparison of the cost of performing a resultant computation guided by the decompositions shown in \cref{fig:doubleBananaW4}.}\label{tbl:dblBan}
\end{table}

\begin{figure}[h]
	\centering
	\begin{tabular}{ccc}
		\includegraphics[width=.35\textwidth]{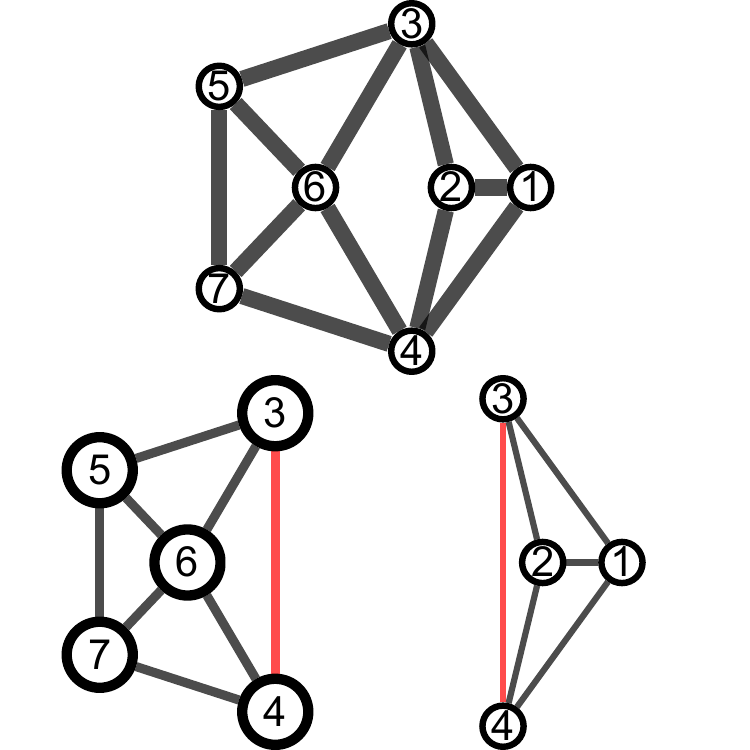} & \hspace{.1cm} &
		\includegraphics[width=.35\textwidth]{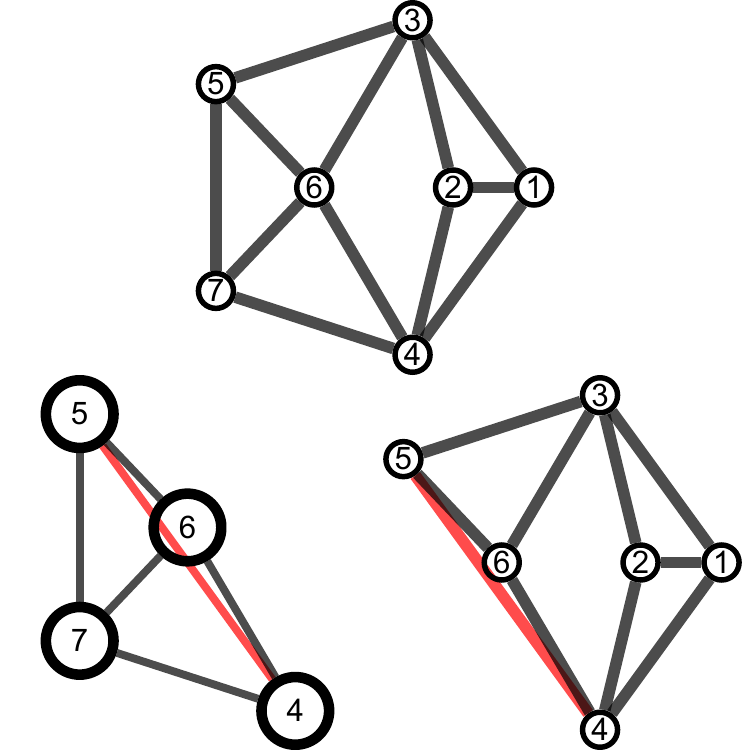}
	\end{tabular}
	\caption{Two \crds of a 2-connected circuit on 7 vertices. The elimination edge is shown in red. Left: a 2-split. Right: a \crd with the triangle 456 in common.}\label{fig:more2conn}
\end{figure}

\newpage
\noindent\textbf{Example 2}. \cref{tbl:tutte} compares the computational time of the circuit polynomial for the 2-connected circuit on 7 vertices shown in \cref{fig:more2conn}. On the left of the figure a 2-split is show, and on the right the two circuits in the \crd have a triangle in common.

The 2-split was computed in 16.94 seconds, compared to 23.5 seconds it took to compute the \crd with the triangle in common.

\begin{table}[h]
	\centering
	\begin{tabular}[t]{|l|c|}
		\hline
		\multicolumn{2}{|c|}{Elimination strategies guided by \cref{fig:more2conn}}\\
		\hline
		Syl Det Dimension for both & $6\times 6$\\
		\hline
		Syl Resultant for both & Hom.\ poly.\ of hom.\ deg.\ 20 \\
		\hline
		\# of mon.\ terms & 1 053 933 \\
		\hline
		Computational time for the 2-split & 16.94 seconds \\
		\hline
		Computational time for the \crd & \multirow{2}{*}{23.50 seconds}\\
		with a common triangle & \\
		\hline
	\end{tabular}
	\caption{Comparison of the cost of performing resultant computations guided \cref{fig:more2conn}. A 2-split and a \crd with a triangle in common are compared.}\label{tbl:tutte}
\end{table}

\begin{figure}[h]
	\centering
	\begin{tabular}{ccc}
		\includegraphics[width=.35\textwidth]{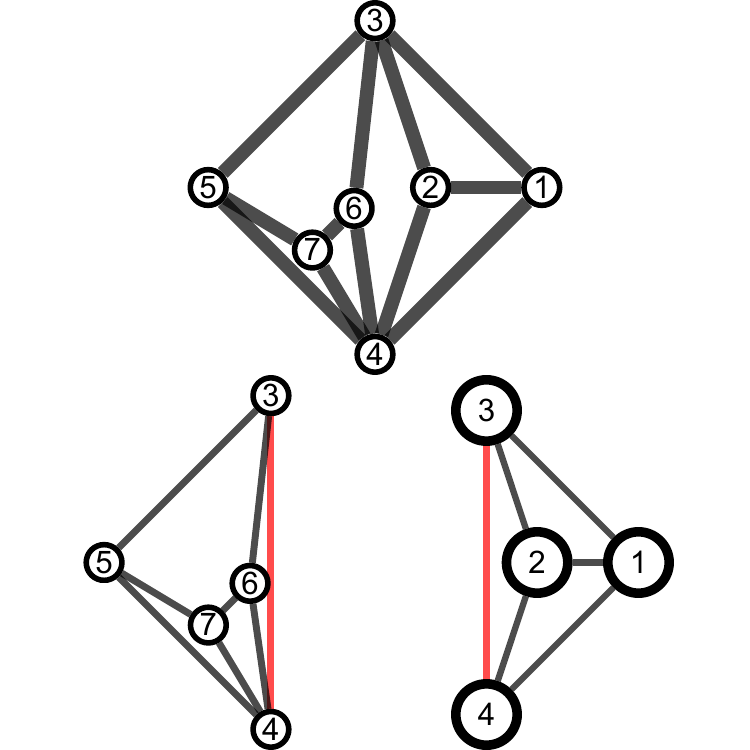} & \hspace{.1cm} & \includegraphics[width=.35\textwidth]{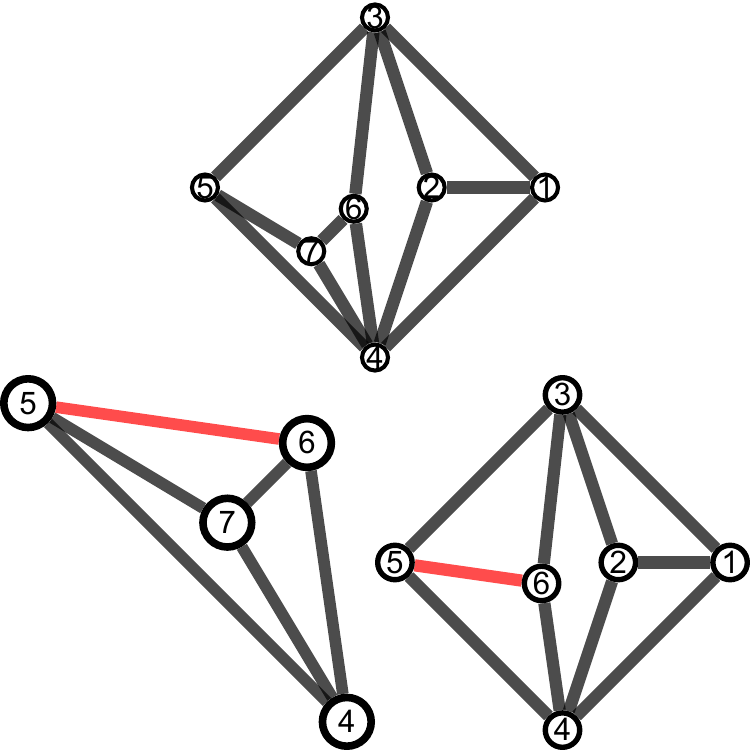}
	\end{tabular}
	\caption{Two \crds of a 2-connected circuit on 7 vertices. The elimination edge is shown in red. Left: a 2-split. Right: a \crd with the triangle 456 in common.}\label{fig:more2connCont}
\end{figure}

\newpage
\noindent\textbf{Example 3}. \cref{tbl:tutte2} compares the computational time of the circuit polynomial for the 2-connected circuit on 7 vertices shown in \cref{fig:more2connCont}. On the left of the figure a 2-split is shown, and on the right the two circuits in the \crd have a triangle in common.

The 2-split was computed in 42.19 seconds, compared to 59.88 seconds it took to compute the \crd with the triangle in common.

\begin{table}[h]
	\centering
	\begin{tabular}[t]{|l|c|}
		\hline
		\multicolumn{2}{|c|}{Elimination guided by \cref{fig:more2connCont}}\\
		\hline
		Syl Det Dimension for both & $6\times 6$\\
		\hline
		Syl Resultant for both & Hom.\ poly.\ of hom.\ deg.\ 20 \\
		\hline
		\# of mon.\ terms & 2 579 050 \\
		\hline
		Computational time for the 2-split & 42.19 seconds \\
		\hline
		Computational time for the \crd & \multirow{2}{*}{59.88 seconds}\\
		with a common triangle & \\
		\hline
	\end{tabular}\vspace{.5cm}
	\caption{Comparison of the cost of performing resultant computations guided \cref{fig:more2connCont}. A 2-split and a \crd with a triangle in common are compared.}\label{tbl:tutte2}
\end{table}

\begin{figure}[h]
	\centering
	\begin{tabular}{ccc}
		\includegraphics[width=.35\textwidth]{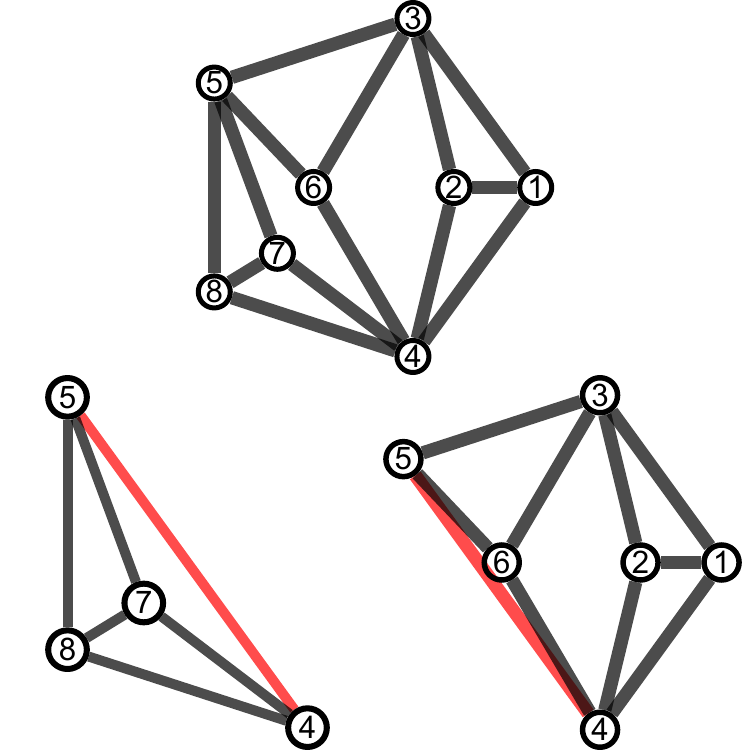} & \hspace{.1cm} & \includegraphics[width=.35\textwidth]{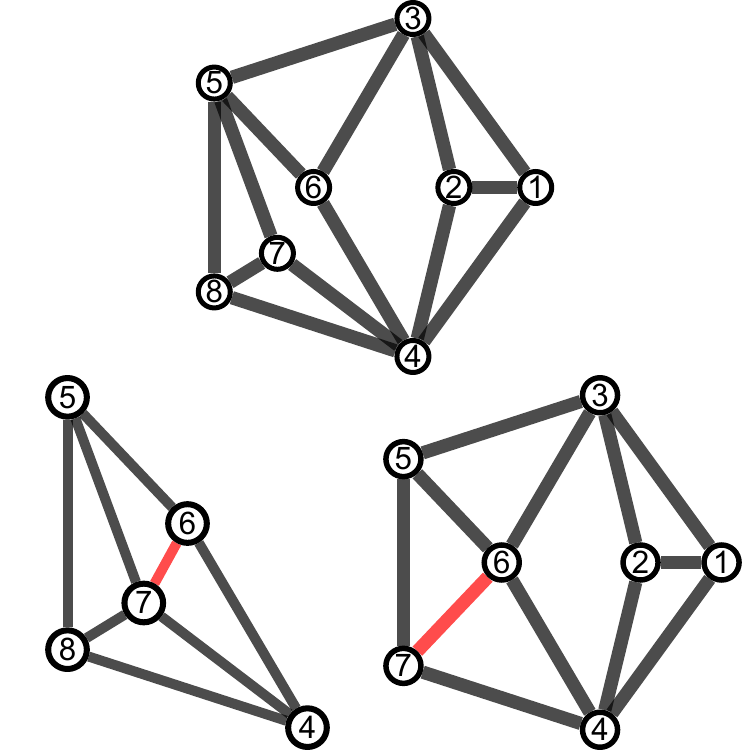}
	\end{tabular}
	\caption{Two \crds of a 2-connected circuit on 8 vertices. The elimination edge is shown in red. Left: a 2-split. Right: a \crd with the double triangle 4675 in common.}\label{fig:8v1}
\end{figure}

\newpage
\noindent\textbf{Example 4}. \cref{tbl:tutte3} compares the computational time of the circuit polynomial for the 2-connected circuit on 8 vertices shown in \cref{fig:8v1}. On the left of the figure a 2-split is shown, and on the right the two circuits in the \crd have a double triangle in common.

The 2-split was computed in 76.67 seconds, while the computation guided by the \crd on the right crashed.

\begin{table}[h]
	\centering
	\begin{tabular}[t]{|l|c|}
		\hline
		\multicolumn{2}{|c|}{Elimination guided by the 2-split in \cref{fig:8v1}}\\
		\hline
		Syl Det Dimension& $6\times 6$\\
		\hline
		Syl Resultant& Hom.\ poly.\ of hom.\ deg.\ 20 \\
		\hline
		\# of mon.\ terms & 3 413 204 \\
		\hline
		Computational time  & 76.67 seconds \\
		\hline
	\end{tabular}\vspace{.5cm}
	\begin{tabular}[t]{|l|c|}
		\hline
		\multicolumn{2}{|c|}{
			Elimination guided by the \crd 
		}\\
		\multicolumn{2}{|c|}{
			with a double triangle in common in \cref{fig:8v1}
		}\\
		\hline
		Syl Det Dimension for both & $12\times 12$\\
		\hline
		Syl Resultant & Hom.\ poly.\ of hom.\ deg.\ 112 \\
		\hline
		\# of mon.\ terms & Crashed \\
		\hline
		Computational time  & Crashed \\
		\hline
	\end{tabular}
	\caption{Comparison of the cost of performing resultant computations guided \cref{fig:8v1}. A 2-split and a \crd with a double triangle in common are compared.}\label{tbl:tutte3}
\end{table}

\begin{figure}[h]
	\centering
		\begin{tabular}{ccc}
				\includegraphics[width=.35\textwidth]{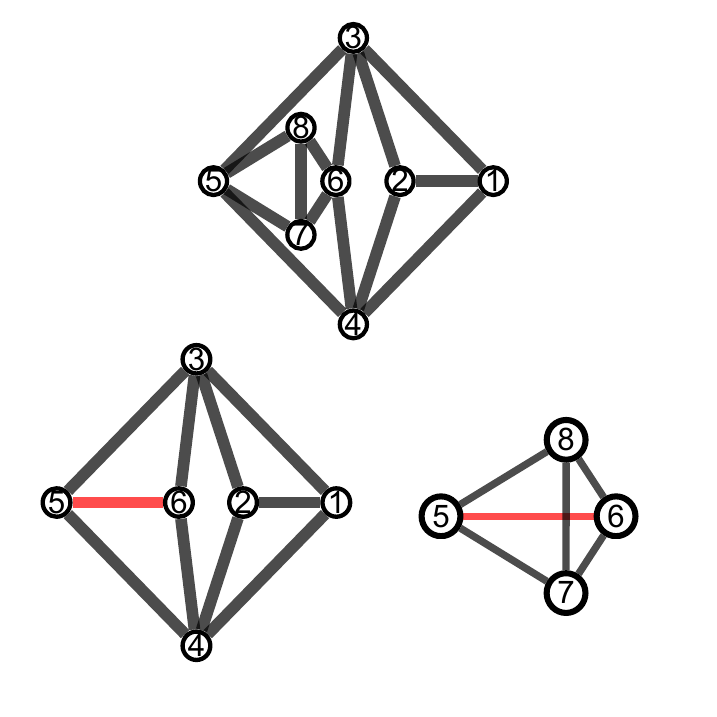} & \hspace{.1cm} & \includegraphics[width=.35\textwidth]{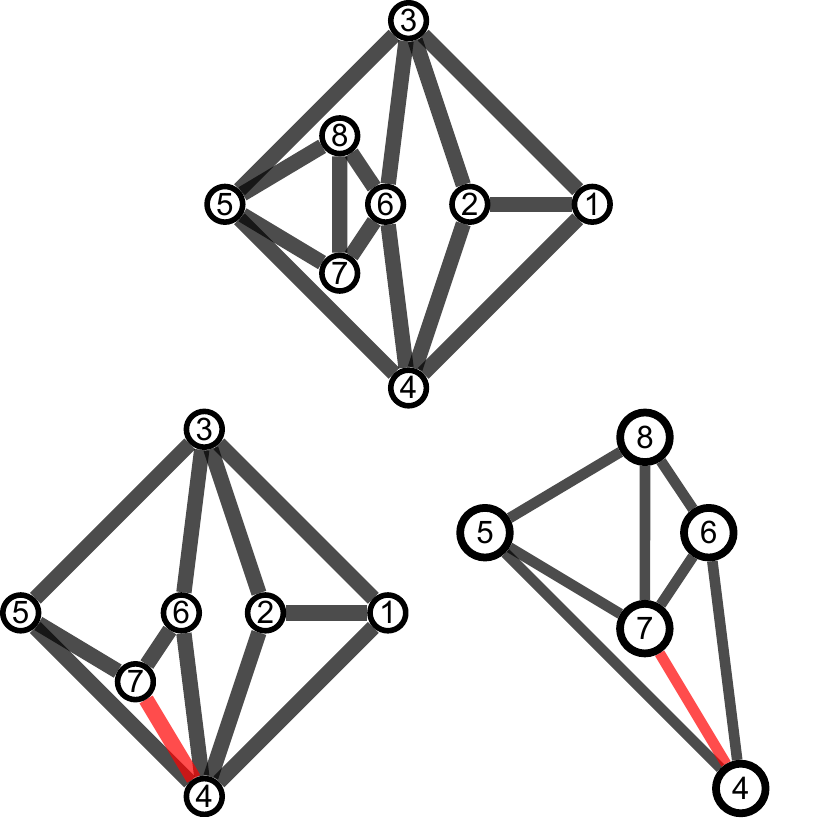}
		\end{tabular}
	\caption{Two \crds of a 2-connected circuit on 8 vertices. The elimination edge is shown in red. Left: a 2-split. Right: a \crd with the double triangle 4675 in common.}\label{fig:8v2}
\end{figure}

\newpage
\noindent\textbf{Example 5}. \cref{tbl:tutte4} compares the computational time of the circuit polynomial for the 2-connected circuit on 8 vertices shown in \cref{fig:8v2}. On the left of the figure a 2-split is shown, and on the right the two circuits in the \crd have a double triangle in common.

The 2-split was computed in 204.71 seconds, while the computation guided by the \crd on the right crashed. 

\begin{table}[h]
	\centering
	\begin{tabular}[t]{|l|c|}
		\hline
		\multicolumn{2}{|c|}{Elimination guided by the 2-split in \cref{fig:8v2}}\\
		\hline
		Syl Det Dimension& $6\times 6$\\
		\hline
		Syl Resultant& Hom.\ poly.\ of hom.\ deg.\ 20 \\
		\hline
		\# of mon.\ terms & 9 223 437 \\
		\hline
		Computational time  & 204.71 seconds \\
		\hline
	\end{tabular}\vspace{.5cm}
	\begin{tabular}[t]{|l|c|}
		\hline
		\multicolumn{2}{|c|}{
			Elimination guided by the \crd 
			}\\
		\multicolumn{2}{|c|}{
				with a double triangle in common in \cref{fig:8v2}
			}\\
		\hline
		Syl Det Dimension for both & $12\times 12$\\
		\hline
		Syl Resultant & Hom.\ poly.\ of hom.\ deg.\ 112 \\
		\hline
		\# of mon.\ terms & Crashed \\
		\hline
		Computational time  & Crashed \\
		\hline
	\end{tabular}
	\caption{Comparison of the cost of performing resultant computations guided \cref{fig:8v2}. A 2-split and a \crd with a double triangle in common are compared.}\label{tbl:tutte4}
\end{table}

\section{Concluding remarks}\label{sec:concl}
 Evidence presented in Examples 1--5 shows that 2-splits outperform other \crds in those cases and supports our heuristic for choosing a good elimination strategy for computing the circuit polynomials of 2-connected circuits. Whether this extends to larger examples remains an open problem to be tackled by a comprehensive project in which new data structures and algorithms optimized for very large polynomials would have to be considered.

\newpage
\bibliography{references}

\end{document}